# A Computational Procedure for solving a Non-Convex Multi-Objective Quadratic Programming under Fuzzy Environment

**Shashi Aggarwal***
Department of Mathematics
Miranda House, University of Delhi, Delhi-110007, India
shashi60@gmail.com

**Uday Sharma**
Department of Mathematics
University of Delhi, Delhi-110007, India
udaysharm88@yahoo.com

The purpose of this paper is to study a non-convex fuzzy multi-objective quadratic programming problem, in which both the technological coefficients and resources are fuzzy with nonlinear membership function. A computational procedure to find a fuzzy efficient solution of this problem is developed. A numerical example is given to illustrate the procedure.

## 1. Introduction

This paper studies the problem of maximizing a number of quadratic objective functions subject to linear and bound constraints under fuzzy environment:

**(NFMOQPP)**

$$\text{Max } Z_1(x) = c_1^t x + \frac{1}{2} x^t Q_1 x$$

$$\text{Max } Z_2(x) = c_2^t x + \frac{1}{2} x^t Q_2 x$$

$$\vdots$$

$$\text{Max } Z_k(x) = c_k^t x + \frac{1}{2} x^t Q_k x$$

$$\text{Subject to } \quad x \in X = \{x \in R^n : \tilde{A} x \le \tilde{b},\ \tilde{l} \le x \le \tilde{u},\ x \ge 0\},$$

where $c_q$, q = 1,2, …k are n-dimensional cost vectors; $Q_q$, q = 1, 2…k are $n \times n$ symmetric positive semi-definite matrices; $\tilde{A}$ is the $m \times n$ constraint fuzzy matrix of technological coefficients; $\tilde{b}$ is the m-dimensional fuzzy resource vector, n-dimensional fuzzy vectors $\tilde{l}$ and $\tilde{u}$ are lower and upper bounds respectively on n-dimensional decision vector x.

* Corresponding Author.

**Remark 1.1.** The problem **(NFMOQPP)** is a non-convex programming problem because $Q_q$ , q = 1, 2…k are assumed to be $n \times n$ symmetric positive semi-definite matrices and our problem is maximization.

In the literature section, a greater part of the study focused on convex programming problem; as in convex programming local optima is global optima but in case of non-convex programming local optima may not be global optima. Orden (1963) was the first one who considered the maximization of a non-concave quadratic programming. He derived the necessary and sufficient conditions for the maximization of non-concave quadratic programming. Ritter (1966) extended Orden's work and gave a different approach to find the maximization of non-concave quadratic programming. Murty and Kabadi (1987) discussed that non-convex problem is an N-P Hard problem and it is very hard to get the global maxima. They have conferred that descent algorithms are quite practical algorithms for dealing with non-convex nonlinear problem. Hanson (1990) used generalized invexity to derive the necessary and sufficient conditions in non-convex quadratic programming for global minima. Ye (1992) applied the affine scaling algorithm to find the optimal solution for the non-convex quadratic programming. Burer and Vandenbussche (2008) has developed a finite branch and bound algorithm for non-convex quadratic programming via semi definite relaxation. Chen and Burer (2011) extended the work of Burer and Vandenbussche to get the global optima of non-convex problem with linear and bound constraints.

The purpose of this paper is to develop a computational procedure for more general non-convex quadratic problem similar to that of Chen and Burer (2011) and Orden (1963). They have focused on the single objective non-convex quadratic programming problem. But here in this paper, we have concentrated on the problem **(NFMOQPP)** in which instead of single objective we have taken multiple objectives which are fuzzified during the procedure in order to obtain the bounds on them; the crisp linear and bound constraints have been replaced by fuzzy linear and bound constraints.

In real world, the concept of decision-making takes place in an environment in which the objectives and constraints are not known precisely. Such situations can be tackled efficiently with the help of fuzzy set theory. Fuzzy sets were first introduced by Zadeh (1965). These sets were used to introduce the concept of decision-making in a fuzzy environment by Bellman and Zadeh (1970). They have defined appropriate aggregation of fuzzy sets for the fuzzy decision. Zimmermann (1978) has used fuzzy decision concept in the fuzzy linear programming for several objectives. Guu and Wu (1997) extended the Zimmermann's approach to two-phase approach for solving the multi-objective linear programming in the fuzzy environment. Several authors have studied linear programming in the fuzzy environment and applied to real world problems like transportation, production planning etc. A fuzzy decision is based on the intersection of membership functions of the goals and constraints. Most of the authors have used linear membership functions as they are easy to tackle. Leberling (1981) has defined the hyperbolic membership function. Li and Lee (1991) have defined the exponential membership function and Yang, Ignizio and Kim (1991) have defined piecewise nonlinear membership function. In our paper, we have used trigonometric membership function in terms of sin, for objective functions and some of the constraints and used LINGO 9.0 to find the fuzzy efficient solution of **(NFMOQPP)**.

The paper is organized as follows: Section 2 consists of preliminaries, which contain a useful definition and some basic membership functions of fuzzy technological coefficients and fuzzy resource vector; fuzzy efficient solution and nonlinear membership functions are described in Section 3; in Section 4, we have discussed the solution methodology for solving non-convex fuzzy multi-objective quadratic programming problem; and illustrative numerical example to explain the solution methodology is given in Section 5.

## 2. Preliminaries

In this section, we give a basic definition and membership functions corresponding to **(NFMOQPP)**.

General form of Non-convex Multi-objective Quadratic Programming problem is given by

**(NMOQPP)** $\quad \text{Max } Z_1(x) = c_1^t x + \frac{1}{2} x^t Q_1 x$

$$\text{Max } Z_2(x) = c_2^t x + \frac{1}{2} x^t Q_2 x$$

$$\vdots$$

$$\text{Max } Z_k(x) = c_k^t x + \frac{1}{2} x^t Q_k x$$

$$\text{Subject to} \quad x \in X = \{x \in R^n : Ax \le b, x \ge 0\},$$

where $c_q$, q = 1,2, …k are n-dimensional cost vectors; $Q_q$, q = 1, 2… k are $n \times n$ symmetric positive semi-definite matrices; A is the $m \times n$ constraint matrix of technological coefficients; b is the m-dimensional resource vector and x is n-dimensional decision vector.

**Definition 2.1.** A point $x^* \in X$ is said to be Pareto optimal solution of **(NMOQPP)** if there does not exist any $x \in X$ such that

$$Z_q(x) \ge Z_q(x^*) \ \forall \ q = 1, 2\ldots k \quad \text{and} \quad Z_r(x) > Z_r(x^*) \text{ for at least one } r = 1, 2\ldots k. \tag{1}$$

In the problem **(NFMOQPP)**, $\tilde{A}$ is the $m \times n$ constraint fuzzy matrix of technological coefficients; $\tilde{b}$ is the m-dimensional fuzzy resource vector, $\tilde{l}$ and $\tilde{u}$ are n-dimensional fuzzy vectors and x is n-dimensional decision vector. Their membership functions are given as below:

1. The membership function of $\tilde{b}$ :

$$\mu_{\tilde{b}}(y) = \begin{cases} 1, & y \le b_i \\ \dfrac{b_i + p_i - y}{p_i}, & b_i \le y \le b_i + p_i \\ 0, & y \ge b_i + p_i \end{cases} \tag{2}$$

where $y \in R$ and $p_i > 0$ is the tolerance level of $b_i$ for all i=1, 2… m.

2. The membership function of the fuzzy matrix $\tilde{A}$ :

$$\mu_{\tilde{A}}(y) = \begin{cases} 1, & y \le a_{ij} \\ \dfrac{a_{ij} + d_{ij} - y}{d_{ij}}, & a_{ij} \le y \le a_{ij} + d_{ij} \\ 0, & y \ge a_{ij} + d_{ij} \end{cases} \tag{3}$$

where $y \in R$ and $d_{ij} > 0$ is the tolerance level of $a_{ij}$ for i=1,2…….,m and j=1,2,……..,n.

3. The membership function of $\tilde{u}$ :

$$\mu_{\tilde{u}}(y) = \begin{cases} 1, & y \le u_j \\ \dfrac{u_j + t_j - y}{t_j}, & u_j \le y \le u_j + t_j \\ 0, & y \ge u_j + t_j \end{cases} \tag{4}$$

where $y \in R$ and $t_j > 0$ is the tolerance level of $u_j$ for all j=1, 2,…n.

4. The membership function of $\tilde{l}$ :

$$\mu_{\tilde{l}}(y) = \begin{cases} 1, & l_j \le y \\ \dfrac{y - l_j + r_j}{r_j}, & l_j - r_j \le y \le l_j \\ 0, & y \le l_j - r_j \end{cases} \tag{5}$$

where $y \in R$ and $r_j > 0$ is the tolerance level of $l_j$ for all j=1, 2…, n.

## 3. Fuzzy Efficient Solution

Werners (1987) has given a definition of fuzzy efficient solution for fuzzy multi-objective linear programming problem. As in our problem **(NFMOQPP)**, the constraint set X is described exclusively by fuzzy constraints, the varying degree of feasibility should be taken into account by consideration of efficient solutions, for additional dependencies between the individual goal values and the degree of membership to the region of feasible solutions can arise. That means we want to emphasize not only the achievement of maximum value of objective functions, but also the highest membership grade of fuzzy constraints in a flexible region. In order to take care of fuzzy concept, we will extend the definition of Werners to **(NFMOQPP)**.

Let the constraint sets of **(NFMOQPP)** be denoted by

$C_i(x) = \{ x \in R^n : \tilde{A}_i x \le \tilde{b}_i \}$, i=1, 2……., m

$B_j(x) = \{ x \in R^n : x_j \le \tilde{u}_j \}, \quad j = 1, 2, \ldots, n$

$B'_j(x) = \{ x \in R^n : x_j \ge \tilde{l}_j \}, \quad j = 1, 2, \ldots, n$

Define the trigonometric membership function of i $^{th}$ constraint $C_i(x)$ as:

$$\mu\, C_i(x) = \begin{cases} 0, & b_i \le \sum_{j=1}^{n} a_{ij} x_j \\ \sin\left[\left(\dfrac{b_i - \sum_{j=1}^{n} a_{ij} x_j}{\sum_{j=1}^{n} d_{ij} x_j + p_i}\right)\dfrac{\pi}{2}\right], & \sum_{j=1}^{n} a_{ij} x_j \le b_i \le \sum_{j=1}^{n} (a_{ij} + d_{ij}) x_j + p_i \\ 1, & b_i \ge \sum_{j=1}^{n} (a_{ij} + d_{ij}) x_j + p_i \end{cases} \tag{6}$$

where $x \in R^n$ and $d_{ij} > 0$ and $p_i > 0$ are respectively the tolerance level of the technological coefficients and resources for i=1, 2…,m and j= 1, 2…,n.

This trigonometric membership function has the following properties:

1. $\mu(C_i(\mathrm{x}))$ is a nonlinear and monotonically decreasing function .
2. $0 \le \mu(C_i(\mathrm{x})) \le 1 \ \forall\, \mathrm{x} \in \mathrm{R}^n$.
3. $\mu(C_i(\mathrm{x}))$ is a concave function on the set $\left\{ \mathrm{x} \in R^n : \sum_{j=1}^{n} a_{ij}\mathrm{x}_j \le \mathrm{b}_i \right\}$.

We define the linear membership function of the constraint $B_j(\mathrm{x})$ , j=1, 2…n as:

$$\mu(B_j(\mathrm{x})) = \begin{cases} 1 & \mathrm{x}_j \le \mathrm{u}_j \\ \dfrac{u_j + t_j - \mathrm{x}_j}{t_j} & \mathrm{u}_j \le \mathrm{x}_j \le \mathrm{u}_j + \mathrm{t}_j \\ 0 & \mathrm{x}_j \ge \mathrm{u}_j + \mathrm{t}_j \end{cases} \tag{7}$$

where $\mathrm{x} \in R^n$ and $t_j$ >0 is the tolerance level of $u_j$ .

Define the linear membership function of the constraint $B'_j(\mathrm{x})$, j = 1,2,…n as :

$$\mu(B'_j(\mathrm{x})) = \begin{cases} 1 & \mathrm{l}_j \le \mathrm{x}_j \\ \dfrac{\mathrm{x}_j - \mathrm{l}_j + \mathrm{r}_j}{r_j} & \mathrm{l}_j - \mathrm{r}_j \le \mathrm{x}_j \le \mathrm{l}_j \\ 0 & \mathrm{l}_j - \mathrm{r}_j \ge \mathrm{x}_j \end{cases} \tag{8}$$

where $\mathrm{x} \in R^n$ and $r_j$ >0 is the tolerance level of $l_j$.

**Definition 3.1.** The point $\mathrm{x}^* \in \mathrm{X}$ is said to be fuzzy efficient solution for the **(NFMOQPP)** if there does not exist any $\mathrm{x} \in \mathrm{X}$ such that

$$\begin{bmatrix} Z_q(\mathrm{x}) \ge Z_q(\mathrm{x}^*)\ \forall\ \mathrm{q} = 1,2,\ldots\mathrm{k} \ \textbf{and}\ \mathrm{C}_i(\mathrm{x}) \ge \mathrm{C}_i(\mathrm{x}^*)\ \forall\ \mathrm{i} = 1,2,\ldots\mathrm{m}\ \textbf{and} \\ \mathrm{B}_j(\mathrm{x}) \ge \mathrm{B}_j(\mathrm{x}^*)\ \forall\ \mathrm{j} = 1,2,\ldots\mathrm{n}\ \textbf{and}\ \mathrm{B}'_j(\mathrm{x}) \ge \mathrm{B}'_j(\mathrm{x}^*)\ \forall\ \mathrm{j} = 1,2,\ldots\mathrm{n} \end{bmatrix}$$

**and**

$$\begin{bmatrix} Z_q(\mathrm{x}) > Z_q(\mathrm{x}^*) \text{ for at least one q} = 1,2,\ldots\mathrm{k}\ \textbf{or}\ \mathrm{C}_i(\mathrm{x}) > \mathrm{C}_i(\mathrm{x}^*) \text{ for at least one i} = 1,2,\ldots\mathrm{m}\ \textbf{or} \\ \mathrm{B}_j(\mathrm{x}) > \mathrm{B}_j(\mathrm{x}^*) \text{ for at least one j} = 1,2,\ldots\mathrm{n}\ \textbf{or}\ \mathrm{B}'_j(\mathrm{x}) > \mathrm{B}'_j(\mathrm{x}^*) \text{ for at least one j} = 1,2,\ldots\mathrm{n} \end{bmatrix}$$

As mentioned in Werners (1987), similarly here, this definition takes into account that for a fuzzy efficient solution an improvement concerning an objective function can only be reached either at the expense of an additional objective function or at the expense of the membership into the constraints. In fact, we can easily see that the fuzzy efficiency defined above comprises the classical efficiency as special case if each of the $\mu(C_i(\mathrm{x}))$ , $\mu(B_j(\mathrm{x}))$ and $\mu(B'_j(\mathrm{x}))$ is 1 for each $\mathrm{x} \in \mathrm{C}_i(\mathrm{x})$, $\mathrm{x} \in B_j(\mathrm{x})$ and $\mathrm{x} \in B'_j(\mathrm{x})$ respectively.

## 4. Procedure to find a Fuzzy Efficient Solution

Our objective functions are crisp and constraints are fuzzy in nature in **(NFMOQPP)**. To defuzzificate the problem, we will first fuzzify the objective functions, which can be done by solving the following four quadratic programming problems

**(NQPP $^1_q$)** $Z^1_q$ = Max $Z_q(x) = c^t_q x + \frac{1}{2} x^t Q_q x$, q = 1, 2…, k

Subject to $\sum_{j=1}^{n} a_{ij} x_j \leq b_i$, i = 1, 2…, m

$x_j \leq u_j$, j= 1, 2…. n

$l_j \leq x_j$, j= 1, 2…..n.

$x \geq 0$.

**(NQPP $^2_q$)** $Z^2_q$ = Max $Z_q(x) = c^t_q x + \frac{1}{2} x^t Q_q x$, q = 1, 2…, k

Subject to $\sum_{j=1}^{n} (a_{ij} + d_{ij}) x_j \leq b_i$, i = 1, 2…, m

$x_j \leq u_j$, j= 1, 2…. n

$l_j - r_j \leq x_j$, j= 1, 2…..n.

$x \geq 0$.

**(NQPP $^3_q$)** $Z^3_q$ = Max $Z_q(x) = c^t_q x + \frac{1}{2} x^t Q_q x$, q = 1, 2…, k

Subject to $\sum_{j=1}^{n} (a_{ij} + d_{ij}) x_j \leq b_i + p_i$, i = 1, 2…, m

$x_j \leq u_j + t_j$, j= 1, 2…. n

$l_j - r_j \leq x_j$, j= 1, 2…..n.

$x \geq 0$.

**(NQPP $^4_q$)** $Z^4_q$ = Max $Z_q(x) = c^t_q x + \frac{1}{2} x^t Q_q x$, q = 1, 2…, k

Subject to $\sum_{j=1}^{n} a_{ij} x_j \leq b_i + p_i$, i = 1, 2…, m

$x_j \leq u_j + t_j$, j= 1, 2…. n

$l_j \leq x_j$, j= 1, 2…..n.

$x \geq 0$.

We have now four non-convex quadratic problems corresponding to $q^{th}$ (q = 1,2…..,k) objective function. They can be solved by using LINGO 9.0 to get the aspiration level for the $q^{th}$ objective function.

Now we take the best and worst value of the optimal solutions of the four non-convex quadratic programming problems.

Let $Z^L_q$ = Min ($Z^1_q$, $Z^2_q$, $Z^3_q$, $Z^4_q$) and $Z^U_q$ = Max ($Z^1_q$, $Z^2_q$, $Z^3_q$, $Z^4_q$), q=1, 2…, k. (9)

Now taking the interval [$Z^L_q$, $Z^U_q$] as the tolerance level interval for the $q^{th}$ objective (q=1, 2…, k) in **(NFMOQPP)**, we define the trigonometric membership function of the $q^{th}$ objective (q=1, 2…, k) as:

$$\mu\left(Z_q(\text{x})\right) = \begin{cases} 0, & Z_q(\text{x}) \le Z_q^L \\ \sin\left[\left(\dfrac{Z_q(\text{x}) - Z_q^L}{Z_q^U - Z_q^l}\right)\dfrac{\pi}{2}\right], & Z_q^L \le Z_q(\text{x}) \le Z_q^U \\ 1, & Z_q(\text{x}) \ge Z_q^U \end{cases} \tag{10}$$

This trigonometric membership function has the following properties:

1. $\mu(Z_q(\text{x}))$ is a nonlinear and monotonically increasing function .
2. $0 \le \mu(Z_q(\text{x})) \le 1$ in $\left[Z_q^L , Z_q^U\right]$ .
3. $\mu(Z_q(\text{x}))$ is concave function in $\left[Z_q^L , \infty\right)$ .

Now, by using max-min fuzzy decision making approach given by Bellman and Zadeh (1970), we have

$$\mu_D(\text{x}) = \underset{q,i,j}{Min}\left(\mu\left(Z_q(\text{x})\right), \mu\left(C_i(\text{x})\right), \mu(B_j(\text{x})), \mu(B'_j(\text{x}))\right) \tag{11}$$

Then, the optimal fuzzy decision is a solution of the problem

$$\underset{\text{x}\ge 0}{Max}\ \mu_D(\text{x}) = \underset{\text{x}\ge 0}{Max}\ (\underset{q,i,j}{Min}\ (\mu(Z_q(\text{x})), \mu(C_i(\text{x})), \mu(B_j(\text{x})), \mu(B'_j(\text{x})))). \tag{12}$$

The problem (12) is equivalent to the following problem as discussed by Zimmermann (1978)

**(NPP)** Max $\lambda$

$$\begin{aligned} \text{Subject to}\quad & \mu(Z_q(\text{x})) \ge \lambda, \quad q = 1,2....,k \\ & \mu(C_i(\text{x})) \ge \lambda, \quad i = 1,2....,m \\ & \mu(B_j(\text{x})) \ge \lambda, \quad j = 1,2.....n \\ & \mu(B'_j(\text{x})) \ge \lambda, \quad j = 1,2.....n \\ & 0 \le \lambda \le 1, \\ & \text{x} \ge 0. \end{aligned}$$

Now, the above problem is a nonlinear convex programming problem which can be solved by LINGO 9.0. Let ( $\text{x}^*$ , $\lambda^*$ ) be the optimal solution of the problem **(NPP)**.

Now, by two phase method proposed by Guu and Wu (1997), we will solve the following nonlinear problem:

**(NPP1)** $\text{Max} \sum_{q=1}^{k}\lambda_q + \sum_{i=1}^{m}\delta_i + \sum_{j=1}^{n}\theta_j + \sum_{j=1}^{n}\gamma_j$

$$\begin{aligned} \text{Subject to}\quad & \mu\left(Z_q(\text{x})\right) \ge \lambda_q \ge \lambda^*, \quad q=1, 2\ldots, k \\ & \mu\left(C_i(\text{x})\right) \ge \delta_i \ge \lambda^*, \quad i=1,2\ldots, m \\ & \mu(B_j(\text{x})) \ge \theta_j \ge \lambda^*, \quad j=1, 2\ldots.n \\ & \mu(B'_j(\text{x})) \ge \gamma_j \ge \lambda^*, \\ & \lambda^* \le \lambda_q , \delta_i , \theta_j , \gamma_j \le 1, \\ & \text{x} \ge 0. \end{aligned}$$

We solve the above problem by LINGO 9.0. Let ( $\text{x}^\circ, \lambda^\circ, \delta^\circ, \theta^\circ, \gamma^\circ$ ) where $\lambda^\circ = (\lambda_1^\circ, \lambda_2^\circ, .., \lambda_k^\circ)$, $\delta^\circ = (\delta_1^\circ, \delta_2^\circ, ..., \delta_m^\circ)$, $\theta^\circ = \left(\theta_1^\circ, \theta_2^\circ \ldots, \theta_n^\circ\right)$ and $\gamma^\circ = \left(\gamma_1^\circ, \gamma_2^\circ \ldots, \gamma_n^\circ\right)$ be the optimal solution of **(NPP1)**. Now we will prove that this optimal solution ( $\text{x}^\circ, \lambda^\circ, \delta^\circ, \theta^\circ, \gamma^\circ$ ) gives the fuzzy efficient solution of **(NFMOQPP)**.

**Theorem 4.1.** Let ( $x^\circ, \lambda^\circ, \delta^\circ, \theta^\circ, \gamma^\circ$ ) where $\lambda^\circ = (\lambda_1^\circ, \lambda_2^\circ, ........, \lambda_k^\circ)$ , $\delta^\circ = (\delta_1^\circ, \delta_2^\circ, ........, \delta_m^\circ)$ , $\theta^\circ = \theta_1^\circ, \theta_2^\circ \ldots, \theta_n^\circ$ and $\gamma^\circ = \gamma_1^\circ, \gamma_2^\circ \ldots, \gamma_n^\circ$ be the optimal solution of the problem (NPP1), then $x^\circ$ is the fuzzy efficient solution of the problem **(NFMOQPP)**.

**Proof:** Let if possible, $x^\circ$ be not a fuzzy efficient solution of **(NFMOQPP)**, then there exist a $y \in X$ such that

$$\begin{bmatrix} Z_q(y) \geq Z_q(x^\circ)\ \forall\ q = 1,2,...k \textbf{ and } C_i(y) \geq C_i(x^\circ)\ \forall\ i = 1,2,...m \textbf{ and} \\ B_j(y) \geq B_j(x^\circ)\ \forall\ j = 1,2,...n \textbf{ and } B'_j(y) \geq B'_j(x^\circ)\ \forall\ j = 1,2,...n \end{bmatrix}$$

$$\textbf{and} \qquad (13)$$

$$\begin{bmatrix} Z_q(y) > Z_q(x^\circ) \text{ for at least one } q = 1,2,...k \textbf{ or } C_i(y) > C_i(x^\circ) \text{ for at least one } i = 1,2,...m \textbf{ or} \\ B_j(y) > B_j(x^\circ) \text{ for at least one } j = 1,2,...n \textbf{ or } B'_j(y) > B'_j(x^\circ) \text{ for at least one } j = 1,2,...n \end{bmatrix}$$

As $\mu(Z_q(.))$ is the increasing function

$Z_q(y) \geq Z_q(x^\circ) \quad \Rightarrow \mu(Z_q(y)) \geq \mu(Z_q(x^\circ))$ for all q = 1, 2… k.

Since ( $x^\circ, \lambda^\circ, \delta^\circ, \theta^\circ, \gamma^\circ$ ) is the optimal solution and the coefficients in the objective function are positive in problem **(NPP1)**

So, we have

$\lambda_q^\circ = \mu(Z_q(x^\circ))$ q=1, 2…, k

$\delta_i^\circ = \mu(C_i(x^\circ))$ i= 1, 2…, m

$\theta_j^\circ = \mu(B_j(x^\circ))$ j= 1,2…n

$\gamma_j^\circ = \mu(B'_j(x^\circ))$ j= 1,2…n

Now choosing,

$\lambda_q^y = \mu(Z_q(y))$ q=1, 2…, k

$\delta_i^y = \mu(C_i(y))$ i= 1, 2…, m

$\theta_j^y = \mu(B_j(y))$ j= 1,2…n

$\gamma_j^y = \mu(B'_j(y))$ j= 1,2…n

As ( $x^\circ, \lambda^\circ, \delta^\circ, \theta^\circ, \gamma^\circ$ ) is the feasible solution of **(NPP1)** so using inequalities (13) and property of $\mu(Z_q(x))$, we get

$\mu(Z_q(y)) \geq \mu(Z_q(x^\circ)) \geq \lambda_q \geq \lambda^*,$ q=1, 2…, k

$\mu(C_i(y)) \geq \mu(C_i(x^\circ)) \geq \delta_i \geq \lambda^*,$ i=1,2…, m

$\mu(B_j(y)) \geq \mu(B_j(x^\circ)) \geq \theta_j \geq \lambda^*$

$\mu(B'_j(y)) \geq \mu(B'_j(x^\circ)) \geq \gamma_j \geq \lambda^*$ j=1, 2….n

So, we can easily see that (y, $\lambda^y, \delta^y, \theta^y, \gamma^y$ ), where $\lambda^y = (\lambda_1^y, \lambda_2^y, ........, \lambda_k^y)$ , $\delta^y = (\delta_1^y, \delta_2^y, ......, \delta_m^y)$ , $\theta^y = \theta_1^y, \theta_2^y \ldots, \theta_n^y$ and $\gamma^y = \gamma_1^y, \gamma_2^y \ldots, \gamma_n^y$ is the feasible solution of (NPP1).

As (y, $\lambda^y, \delta^y, \theta^y, \gamma^y$ ) and ( $x^\circ, \lambda^\circ, \delta^\circ, \theta^\circ, \lambda^\circ$ ) are solutions of the problem **(NPP1)** and using (13) , we get

$$\sum_{q=1}^{k}\lambda_q^y+\sum_{i=1}^{m}\delta_i^y+\sum_{j=1}^{n}\theta_j^y+\sum_{j=1}^{n}\gamma_j^y=\sum_{q=1}^{k}\mu(Z_q(y))+\sum_{i=1}^{m}\mu(C_i(y))+\sum_{j=1}^{n}\mu\left(B_j(y)\right)+\sum_{j=1}^{n}\mu\left(B'_j(y)\right)$$

$$>\sum_{q=1}^{k}\mu(Z_q(x^\circ))+\sum_{i=1}^{m}\mu(C_i(x^\circ))+\sum_{j=1}^{n}\mu\left(B_j\left(x^\circ\right)\right)+\sum_{j=1}^{n}\mu\left(B'_j\left(x^\circ\right)\right)$$

$$=\sum_{q=1}^{k}\lambda_q^\circ+\sum_{i=1}^{m}\delta_i^\circ+\sum_{j=1}^{n}\theta_j^\circ+\sum_{j=1}^{n}\gamma_j^\circ$$

This is a contradiction as ($x^\circ$, $\lambda^\circ$, $\delta^\circ, \theta^\circ, \lambda^\circ$) is optimal solution of the problem **(NPP1)**.

So $x^\circ$ is the fuzzy efficient solution of **(NFMOQPP)**. ■

### 4.1 Computational Procedure

The ideas discussed above for finding fuzzy efficient solution of **(NFMOQPP)** can be summarized in the form of algorithm as given below.

**Step 1:**

Construct the four problems of the form **(NQPP $^{j}_{q}$)** (j= 1, 2, 3, 4) corresponding to q^th (q= 1, 2… k) objective.

**Step 2:**

Solve them to get the aspiration level for the q $^{th}$ (q= 1, 2, …., k) objective.

**Step 3:**

Determine the aspiration level of the q^th (q = 1, 2… k) objective by

$$Z_q^u = \text{Max}\left(Z_q^1, Z_q^2, Z_q^3, Z_q^4\right) \text{ and } Z_q^l = \text{Min}\left(Z_q^1, Z_q^2, Z_q^3, Z_q^4\right) \qquad \forall\ q = 1,2.....,k.$$

**Step 4:**

Construct the trigonometric membership functions of the form (6) and (10) for the constraints and the objectives respectively. Also linear membership functions of the form (7) and (8) respectively for upper and lower values of x.

**Step 5:**

Construct the problem **(NPP)** and solve it. Let ( $x^*$ , $\lambda^*$ ) be the solution and $\lambda^*$ be the optimal value.

**Step 6:**

Construct the problem **(NPP1)** and solve it. Let ( $x^\circ, \lambda^\circ, \delta^\circ, \theta^\circ, \gamma^\circ$ ) where $\lambda^\circ = (\lambda_1^\circ, \lambda_2^\circ, .., \lambda_k^\circ)$, $\delta^\circ = (\delta_1^\circ, \delta_2^\circ, ..., \delta_m^\circ)$, $\theta^\circ = \left(\theta_1^\circ, \theta_2^\circ \ldots, \theta_n^\circ\right)$ and $\gamma^\circ = \left(\gamma_1^\circ, \gamma_2^\circ \ldots, \gamma_n^\circ\right)$ be the optimal solution. Then $x^\circ$ is the fuzzy efficient solution of problem **(NFMOQPP)**.

## 5. Illustrative Example

**(NFMOQPP)**

$$\text{Max } Z_1(x) = x_1 + 2x_2 + x_1^2 + 2x_2^2$$
$$\text{Max } Z_2(x) = 4x_1 + 7x_2 + 2x_1^2 + 3x_2^2$$

$$\text{Subject to } \tilde{1}x_1 + \tilde{1}x_2 \le \tilde{10},$$
$$\tilde{2}x_1 + \tilde{3}x_2 \le \tilde{25},$$
$$\tilde{2} \le x_1 \le \tilde{9},$$
$$\tilde{2} \le x_2 \le \tilde{8},$$
$$x_1, x_2 \ge 0.$$

Here, $A = (a_{ij}) = \begin{pmatrix} 1 & 1 \\ 2 & 3 \end{pmatrix}$, $(d_{ij}) = \begin{pmatrix} 1 & 1 \\ 1 & 2 \end{pmatrix}$, $(a_{ij} + d_{ij}) = \begin{pmatrix} 2 & 2 \\ 3 & 5 \end{pmatrix}$

$b = \begin{pmatrix} 10 \\ 25 \end{pmatrix}$, $p = \begin{pmatrix} 5 \\ 10 \end{pmatrix}$, $(b+p) = \begin{pmatrix} 15 \\ 35 \end{pmatrix}$, $u = \begin{pmatrix} 9 \\ 8 \end{pmatrix}$, $t = \begin{pmatrix} 3 \\ 2 \end{pmatrix}$, $u + t = \begin{pmatrix} 12 \\ 10 \end{pmatrix}$, $l = \begin{pmatrix} 2 \\ 2 \end{pmatrix}$, $r = \begin{pmatrix} 1 \\ 1 \end{pmatrix}$, $l - r = \begin{pmatrix} 1 \\ 1 \end{pmatrix}$

**Computational Procedure:**

**Step1:**

Construction of the four problems of the form **(NQPP $^{j}_{q}$)** (j= 1, 2, 3, 4) corresponding to q[th] (q= 1, 2) objective as follows:

| | |
|---|---|
| **(NQPP $^{1}_{1}$)** $Z_1^1$ = Max $x_1 + 2x_2 + x_1^2 + 2x_2^2$ Subject to $x_1 + x_2 \le 10$, $2x_1 + 3x_2 \le 25$, $2 \le x_1 \le 9$, $2 \le x_2 \le 8$, $x_1, x_2 \ge 0$. | **(NQPP $^{1}_{2}$)** $Z_2^1$ = Max $4x_1 + 7x_2 + 2x_1^2 + 3x_2^2$ Subject to $x_1 + x_2 \le 10$, $2x_1 + 3x_2 \le 25$, $2 \le x_1 \le 9$, $2 \le x_2 \le 8$, $x_1, x_2 \ge 0$. |
| **(NQPP $^{2}_{1}$)** $Z_1^2$ = Max $x_1 + 2x_2 + x_1^2 + 2x_2^2$ Subject to $2x_1 + 2x_2 \le 10$, $3x_1 + 5x_2 \le 25$, $1 \le x_1 \le 9$, $1 \le x_2 \le 8$, $x_1, x_2 \ge 0$. | **(NQPP $^{2}_{2}$)** $Z_2^2$ = Max $4x_1 + 7x_2 + 2x_1^2 + 3x_2^2$ Subject to $2x_1 + 2x_2 \le 10$, $3x_1 + 5x_2 \le 25$, $1 \le x_1 \le 9$, $1 \le x_2 \le 8$, $x_1, x_2 \ge 0$. |
| **(NQPP $^{3}_{1}$)** $Z_1^3$ = Max $x_1 + 2x_2 + x_1^2 + 2x_2^2$ Subject to $2x_1 + 2x_2 \le 15$, $3x_1 + 5x_2 \le 35$, $1 \le x_1 \le 12$, $1 \le x_2 \le 10$, $x_1, x_2 \ge 0$ | **(NQPP $^{3}_{2}$)** $Z_2^3$ = Max $4x_1 + 7x_2 + 2x_1^2 + 3x_2^2$ Subject to $2x_1 + 2x_2 \le 15$, $3x_1 + 5x_2 \le 35$, $1 \le x_1 \le 12$, $1 \le x_2 \le 10$, $x_1, x_2 \ge 0$ |
| **(NQPP $^{4}_{1}$)** $Z_1^4$ = Max $x_1 + 2x_2 + x_1^2 + 2x_2^2$ Subject to $x_1 + x_2 \le 15$, $2x_1 + 3x_2 \le 35$, $2 \le x_1 \le 12$, $2 \le x_2 \le 10$, $x_1, x_2 \ge 0.$ | **(NQPP $^{4}_{2}$)** $Z_2^4$ = Max $4x_1 + 7x_2 + 2x_1^2 + 3x_2^2$ Subject to $x_1 + x_2 \le 15$, $2x_1 + 3x_2 \le 35$, $2 \le x_1 \le 12$, $2 \le x_2 \le 10$, $x_1, x_2 \ge 0.$ |

**Step 2:**

Now by solving the above quadratic programming problems using LINGO 9.0, we get

$Z_1 = ( Z_1^1 , Z_1^2 , Z_1^3 , Z_1^4 )$ = (118, 42, 96.72 and 228.75)

$Z_2 = ( Z_2^1 , Z_2^2 , Z_2^3 , Z_2^4 )$ = (212, 82, 173.68 and 392.5)

**Step 3:**

Now using (9) for q = 1, 2, we get the interval $[Z_q^L , Z_q^U]$ as the aspiration level interval for the $q^{th}$ objective (q=1, 2) in **(NFMOQPP)** as

$[ Z_1^L , Z_1^U ]$ = [42, 228.75] and $[ Z_2^L , Z_2^U ]$ = [82, 392.5]

**Step 4:**

Construction of the membership functions of the form (6) to (8) and (10) for the constraints and objectives, as follows:

$$\mu(C_1(x)) = \begin{cases} 0, & 10 \le x_1 + x_2 \\ \sin\left[\left(\dfrac{10 - x_1 - x_2}{x_1 + x_2 + 5}\right)\dfrac{\pi}{2}\right], & x_1 + x_2 \le 10 \le 2x_1 + 2x_2 + 5 \\ 1, & 10 \ge 2x_1 + 2x_2 + 5 \end{cases}$$

$$\mu(C_2(x)) = \begin{cases} 0, & 25 \le 2x_1 + 3x_2 \\ \sin\left[\left(\dfrac{25 - 2x_1 - 3x_2}{x_1 + 2x_2 + 10}\right)\dfrac{\pi}{2}\right], & 2x_1 + 3x_2 \le 25 \le 3x_1 + 5x_2 + 10 \\ 1, & 25 \ge 3x_1 + 5x_2 + 10 \end{cases}$$

$$\mu(B_1(x)) = \begin{cases} 1, & x_1 \le 9 \\ \dfrac{12 - x_1}{3}, & 9 \le x_1 \le 12 \\ 0, & x_1 \ge 12 \end{cases}$$

$$\mu(B_2(x)) = \begin{cases} 1, & x_2 \le 8 \\ \dfrac{10 - x_2}{2}, & 8 \le x_2 \le 10 \\ 0, & x_2 \ge 10 \end{cases}$$

$$\mu(B_1'(x)) = \begin{cases} 0, & x_1 \le 1 \\ \dfrac{(x_1 - 1)}{1}, & 1 \le x_1 \le 2 \\ 1, & x_1 \ge 2 \end{cases}$$

$$\mu(B_2'(x)) = \begin{cases} 0, & x_2 \le 1 \\ \dfrac{(x_2 - 1)}{1}, & 1 \le x_2 \le 2 \\ 1, & x_2 \ge 2 \end{cases}$$

$$\mu(Z_1(x)) = \begin{cases} 0, & Z_1(x) \le 42 \\ \sin\left[\left(\dfrac{Z_1(x) - 42}{186.75}\right)\dfrac{\pi}{2}\right], & 42 \le Z_1(x) \le 228.75 \\ 1, & Z_1(x) \ge 228.75 \end{cases}$$

$$\mu(Z_2(x)) = \begin{cases} 0, & Z_2(x) \le 82 \\ \sin\left[\left(\dfrac{Z_2(x) - 82}{310.5}\right)\dfrac{\pi}{2}\right], & 82 \le Z_2(x) \le 392.5 \\ 1, & Z_2(x) \ge 392.5 \end{cases}$$

Figure 1

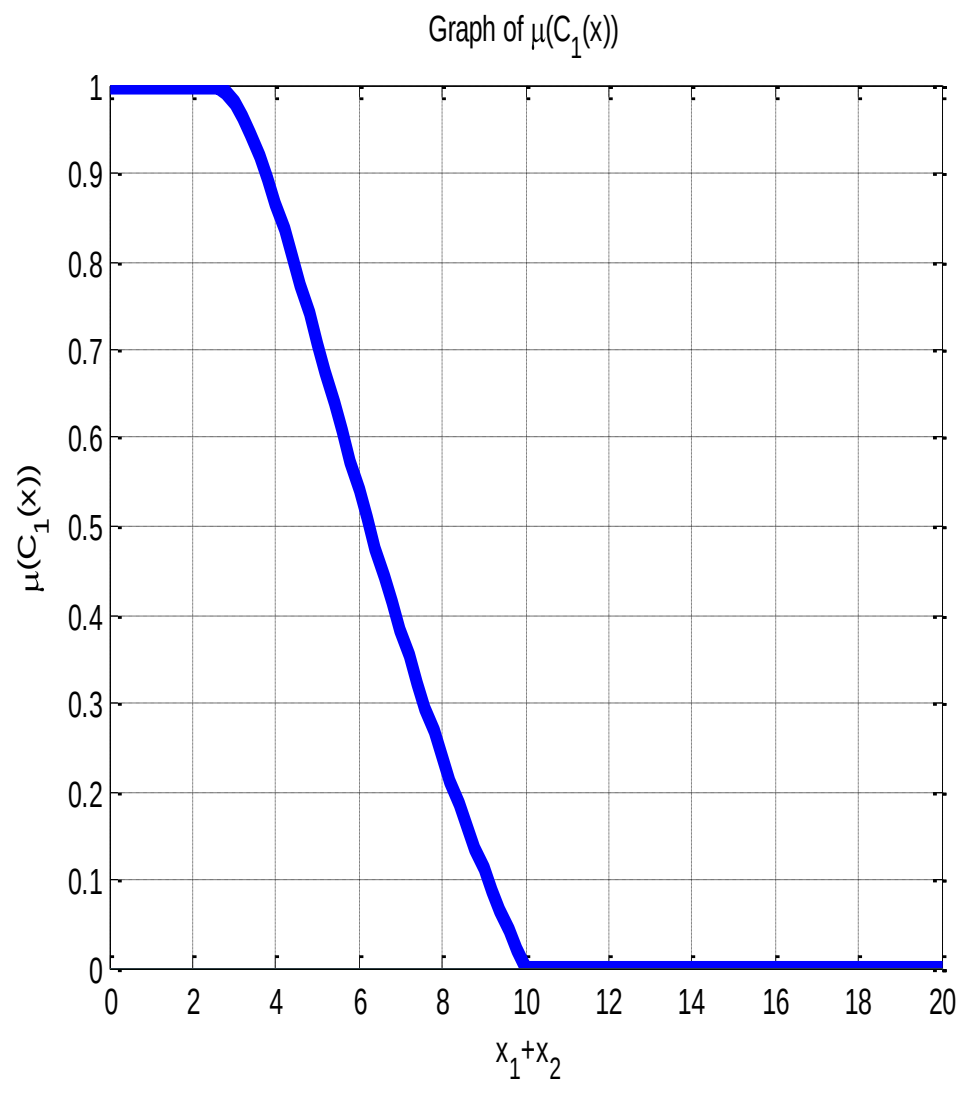

Figure 2

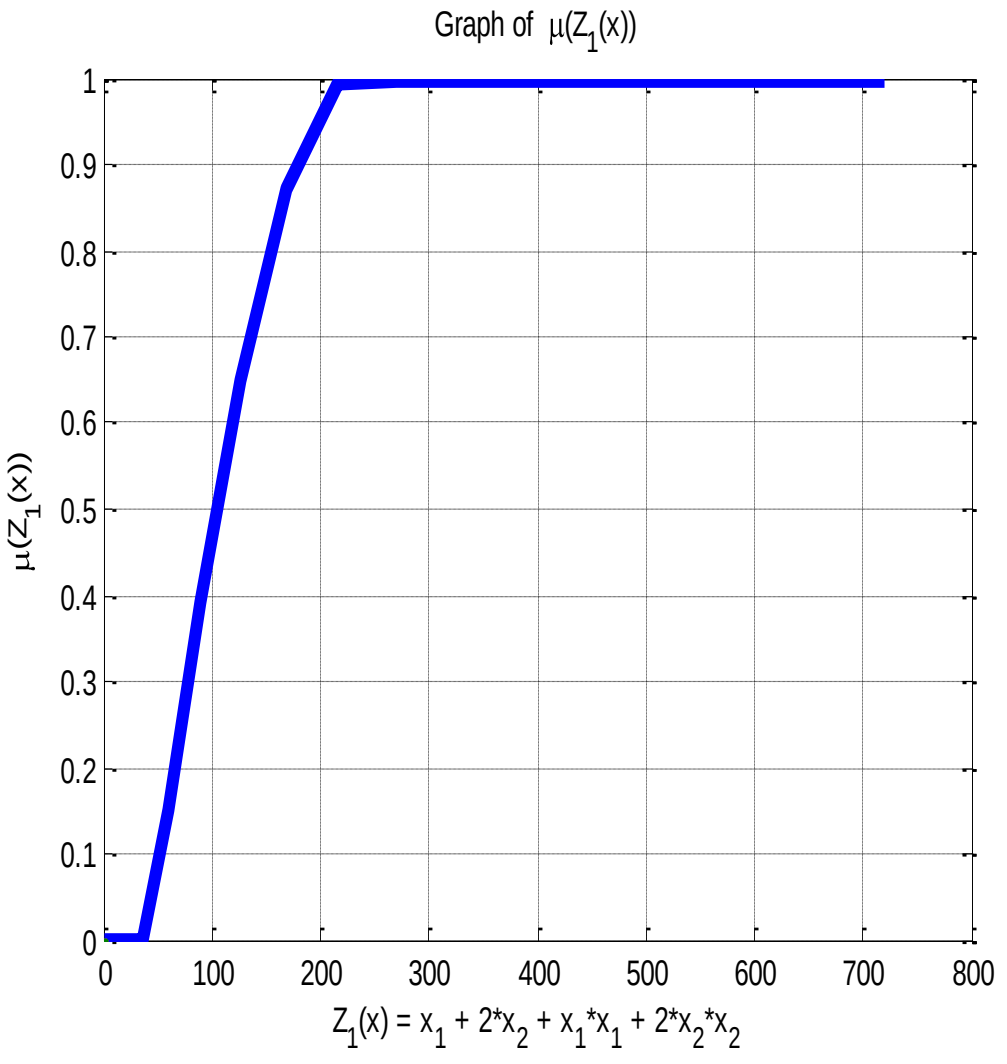

From the graph drawn in figure 1corresponding to the membership function $\mu(C_1(x))$ we can observe that

1. $\mu(C_1(x))$ is a nonlinear and monotonically decreasing function .
2. $0 \le \mu(C_1(x)) \le 1 \;\forall\, x \in R^n$.
3. $\mu(C_1(x))$ is a concave function on the set $\{x \in R^n : x_1 + x_2 \le 10\}$ .

Also from the graph drawn in figure 2 corresponding to the membership function $\mu(Z_1(x))$ we can observe that

1. $\mu(Z_1(x))$ is a nonlinear and monotonic increasing function .
2. $0 \le \mu(Z_1(x)) \le 1 \;\forall\, x \in R^n$ .
3. $\mu(Z_1(x))$ is concave function on the interval $[42, \infty)$ .

Similar properties can be observed for all other membership functions with the help of the graphs.

**Step 5:**

Now the problem **(NPP)** becomes

**(NPP)** Max $\lambda$

Subject to
$$\sin\left[\left(\frac{x_1 + 2x_2 + x_1^2 + 2x_2^2 - 42}{186.75}\right)\frac{\pi}{2}\right] \geq \lambda,$$
$$\sin\left[\left(\frac{4x_1 + 7x_2 + 2x_1^2 + 3x_2^2 - 82}{310.5}\right)\frac{\pi}{2}\right] \geq \lambda,$$
$$\sin\left[\left(\frac{10 - x_1 - x_2}{x_1 + x_2 + 5}\right)\frac{\pi}{2}\right] \geq \lambda,$$
$$\sin\left[\left(\frac{25 - 2x_1 - 3x_2}{x_1 + 2x_2 + 10}\right)\frac{\pi}{2}\right] \geq \lambda,$$
$$\left(\frac{12 - x_1}{3}\right) \geq \lambda,$$
$$\left(\frac{10 - x_2}{2}\right) \geq \lambda,$$
$$(x_1 - 1) \geq \lambda,$$
$$(x_2 - 1) \geq \lambda,$$
$$0 \leq \lambda \leq 1,$$
$$x_1,\ x_2 \geq 0.$$

We solve the above nonlinear problem by LINGO 9.0. We get the optimal solution as ($\lambda^* = 0.3310024$, $x_1^* = 1.331002$, $x_2^* = 5.804374$).

**Step 6:**

Now construct the problem **(NPP1)** as follows:

**(NPP1)** Max $\lambda_1 + \lambda_2 + \delta_1 + \delta_2 + \theta_1 + \theta_2 + \gamma_1 + \gamma_2$

Subject to
$$\sin\left[\left(\frac{x_1 + 2x_2 + x_1^2 + 2x_2^2 - 42}{186.75}\right)\frac{\pi}{2}\right] \geq \lambda_1 \geq 0.3310024,$$
$$\sin\left[\left(\frac{4x_1 + 7x_2 + 2x_1^2 + 3x_2^2 - 82}{310.5}\right)\frac{\pi}{2}\right] \geq \lambda_2 \geq 0.3310024,$$
$$\sin\left[\left(\frac{10 - x_1 - x_2}{x_1 + x_2 + 5}\right)\frac{\pi}{2}\right] \geq \delta_1 \geq 0.3310024,$$
$$\sin\left[\left(\frac{25 - 2x_1 - 3x_2}{x_1 + 2x_2 + 10}\right)\frac{\pi}{2}\right] \geq \delta_2 \geq 0.3310024,$$
$$\left(\frac{12 - x_1}{3}\right) \geq \theta_1 \geq 0.3310024,$$
$$\left(\frac{10 - x_2}{2}\right) \geq \theta_2 \geq 0.3310024,$$
$$(x_1 - 1) \geq \gamma_1 \geq 0.3310024,$$

$$(x_2 - 1) \geq \gamma_2 \geq 0.3310024,$$
$$0.3310024 \leq \lambda_1, \lambda_2, \delta_1, \delta_2, \theta_1, \theta_2, \gamma_1, \gamma_2 \leq 1,$$
$$x_1, x_2 \geq 0.$$

We solve the above problem (NPP1) by LINGO 9.0 and get the optimal values as ( $\lambda_1^\circ$ = 0.3310024, $\lambda_2^\circ$ = 0.3401065, $\delta_1^\circ$ = 0.362496 , $\delta_2^\circ$ = 0.3310024 , $\theta_1^\circ$ = 1 , $\theta_2^\circ$ = 1 , $\gamma_1^\circ$ = 0.3310024 , $\gamma_2^\circ$ = 1 , $x_1^\circ$ = 1.331002 , $x_2^\circ$ = 5.804373 ) . Hence the fuzzy efficient solution of the problem (NFMOQPP) is $x_1^\circ$ = 1.3310024, and $x_2^\circ$ = 5.804373. The corresponding values of the objectives are $Z_1(x^\circ)$ = 82.092806 and $Z_2\ x^\circ$ = 150.569993.

**Acknowledgements**

We are indebted for very helpful discussion to Prof. Davinder Bhatia (Rtd.), Department of Operational Research, Faculty of Mathematical Sciences, University of Delhi, Delhi 110007, India. The second author is thankful to the COUNCIL OF SCIENTIFIC AND INDUSTRIAL RESEARCH (CSIR), NEW DELHI, INDIA, for the financial support (09/045(1153)/2012-EMR-I).

## References

Bector C.R. and S. Chandra (2005). *Fuzzy Mathematical programming and fuzzy matrix games book.* Springer.

Bellmann R.E. and L.A. Zadeh (1970). Decision- making in a fuzzy environment. *Management Science*,17(B),No.4,141-164.

Burer S. and D. Vandenbussche (2008). A finite branch-and-bound algorithm for non-convex quadratic programming via semidefinite relaxations. Mathematical Programming. Ser. A, 113, 259-282.

Chen J. and S. Burer (2011). Globally solving non-convex quadratic programming problems via completely positive programming. Mathematics and Computer Science Division, Argonne National Laboratory, Argonne, Illinois.

Guu,S.M. and Y.K.Wu (1997). Weighted coefficients in two-phase approach for solving the multiple objective programming problems. *Fuzzy Sets and Systems*, 85, 45-48.

Gasimov, R. N.(2002).Augmented Lagrangian Duality and nondifferentiable optimization methods in non-convex programming. Journal of Global Optimization, 24, 187-203.

Hanson M.A. (1990). Non-convex quadratic programming. FUS Technical report number M-810, Department of statistics, Florida State University.

Leberling, H. (1981). On finding compromise solutions in multicriteria problems using the fuzzy min-operator. *Fuzzy sets and systems*, 6, 105-118.

Li R. J. and E. Stanley Lee (1991). An exponential membership function for fuzzy multiple objective linear programming. *Computers and mathematics with applications,* 22, No.-12 , 55-60.

Murty K.G. and S.N. Kabadi (1987). Some NP-complete problems in quadratic and nonlinear programming. Mathematical Programming , 39, 117-129.

Orden A. (1963). Minimization of indefinite quadratic functions with linear constraints. Recent advances in mathematical programming, McGraw Hill.

Ritter K. (1966).A method for solving maximum–problems with a non-concave quadratic objective function. Z. Wahrscheinlichkeitstheorie Verw. Geb. 4, 340-351.

Werners B.(1987) . An Interactive fuzzy programming system. *Fuzzy sets and systems,* 23,131-147.

Werners B. (1987). Interactive multiple objective programming subject to flexible constraints. *European Journal of Operational Research*, 31, 342-349.

Yang,T., J.P. Ignizio and Hyun-Joon Kim (1991). Fuzzy programming with non linear membership functions: piecewise linear approximation. *Fuzzy sets and system* , 41, 39-53.

Ye Yinyu (1992). On affine scaling algorithms for non-convex quadratic programming. Mathematical Programming, 56, 285-300.

Zadeh,L.A. (1965). Fuzzy sets. *Information and Control,* 8, 338-353.

Zimmermann, H.J. (1978). Fuzzy programming and linear programming with several objective functions. *Fuzzy Sets and Systems*, 1, 45-55.

**ShashiAggarwal** is currently working as an Associate Professor (on deputation) in Cluster Innovation Centre, University of Delhi, India. She is the permanent faculty of the Department of Mathematics in Miranda House, University of Delhi, India. She received her doctorate degree in the field of Mathematical Programming from University of Delhi, India, in 1993. She has published her research papers in leading journals like European Journal of Operational Research, Optimization, Asia Pacific Journal of Optimization, Opsearch, Indian Journal of Pure and Applied mathematics, Journal of Information & Optimization Science etc.

**Uday Sharma** received her B.Sc and M.Sc from University of Delhi, India, in 2008 and 2010, respectively. He is currently doing Ph.D under the supervision of first author.